\documentclass[a4paper,12pt]{article}

\usepackage[top=30truemm,bottom=25truemm,left=25truemm,right=25truemm]{geometry}
\usepackage[pdftex]{graphicx}
\usepackage{ulem}
\renewcommand{\title}[1]{
	\vspace*{4mm}
	\begin{center}
	\textbf{\Large #1}
	\end{center}
	\smallskip}

\renewcommand{\author}[1]{
	\vspace*{0mm}
	\begin{center}
	#1
	\end{center}
	\smallskip}
\renewcommand{\abstract}[1]{
	\begin{center}
	\parbox{13cm}{\small {\sc Abstract.} #1}
	\end{center}
	\smallskip}

\renewcommand{\thanks}[1]{{
	\renewcommand{\thefootnote}{\fnsymbol{footnote}}
	\vspace*{-5mm}
	\footnote[0]{#1}}}
\newcommand{\address}[1]{\bigskip{\small\noindent #1 \par}}
\newcommand{\email}[1]{{\small\noindent\textit{Email address}: \texttt{#1} \par}}

\usepackage{titlesec} 
\titleformat*{\section}{\large\bfseries}
\titleformat*{\subsection}{\bfseries}
\titleformat*{\subsubsection}{}
\titlespacing{\section}{0pt}{*3}{*1.5}
\titlespacing{\subsection}{0pt}{*2}{*1}
\titlespacing{\subsubsection}{0pt}{*2}{*1}

 \usepackage{amscd, amsmath, amssymb, amsthm}
 \usepackage[initials,short-journals,non-sorted-cites]{amsrefs}
 \usepackage{tikz}
\usetikzlibrary{
  knots,
  hobby,
  decorations.pathreplacing,
  shapes.geometric,
  calc
}

\usepackage{graphicx}
\usepackage{color}

 \theoremstyle{plain}
 \newtheorem{thm}{Theorem}[section]
 \newtheorem{prop}[thm]{Proposition}
 \newtheorem{cor}[thm]{Corollary}

 \theoremstyle{definition}
 \newtheorem{defn}[thm]{Definition}
 \newtheorem{example}[thm]{Example}
 \newtheorem{rem}[thm]{Remark}

\begin{document}

\tikzset{
  knot diagram/every strand/.append style={
    ultra thick,
    black 
  },
  show curve controls/.style={
    postaction=decorate,
    decoration={show path construction,
      curveto code={
        \draw [blue, dashed]
        (\tikzinputsegmentfirst)--(\tikzinputsegmentsupporta)
        node [at end, draw, solid, red, inner sep=2pt]{};
        \draw [blue, dashed]
        (\tikzinputsegmentsupportb)--(\tikzinputsegmentlast)
        node [at start, draw, solid, red, inner sep=2pt]{}
        node [at end, fill, blue, ellipse, inner sep=2pt]{}
        ;
      }
    }
  },
  show curve endpoints/.style={
    postaction=decorate,
    decoration={show path construction,
      curveto code={
        \node [fill, blue, ellipse, inner sep=2pt] at (\tikzinputsegmentlast) {}
        ;
      }
    }
  }
}
\tikzset{->-/.style 2 args={
    postaction={decorate},
    decoration={markings, mark=at position #1 with {\arrow[thick, #2]{>}}}
    },
    ->-/.default={0.5}{}
}

\tikzset{-<-/.style 2 args={
    postaction={decorate},
    decoration={markings, mark=at position #1 with {\arrow[thick, #2]{<}}}
    },
    -<-/.default={0.5}{}
}
\title{Infinitely many pairs of spatial surfaces}

\author{Katsunori Arai}

 \abstract{
  A multiple group rack (MGR) is an algebraic system 
  which is used to construct invariants of spatial surfaces,
  which are compact surfaces embedded in the $3$-sphere $S^{3}$.
  Seifert surfaces for links are spatial surfaces.
  In this paper,
  we present an infinitely many pairs of Seifert surfaces for each link, where each pair satisfies the following condtions:
  (i) their regular neighborhoods in $S^{3}$ are ambiently isotopic,
  (ii) their Seifert matrices are unimodularly congruent, and
  (iii) the two Seifert surfaces are not ambiently isotopic.
  In order to prove (iii),
  we distinguish the Seifert surfaces using the above invariants.
 }




\section{Introduction}
A \textit{Seifert surface} for a link is a compact connected oriented surface in the $3$-sphere $S^{3}$
whose boundary is ambiently isotopic to the link.
It is well known that every link has a Seifert surface \cite{Frankl-Pontrjagin1930,Seifert1935}.
Two Seifert surfaces for a link are \textit{equivalent} if
they are ambiently isotopic in the exterior of the link.
One of the important problems in knot theory is the classification of Seifert surfaces for a given link up to the equvalence.
This problem has been extensively studied in the case of 
minimal genus or incompressible Seifert surfaces. 
Compact surfaces embedded in the $3$-sphere $S^{3}$ is called \textit{spatial surfaces}.
Two spatial surfaces are said to be \textit{equivalent} if they are ambiently isotopic in $S^{3}$.
Throughout this paper,
we assume that 
(1) a spatial surface is oriented and 
(2) each component of a spatial surface is neither a closed disk nor a surface without boundary.
Under the assumptions,
spatial surfaces are Seifert surfaces for their boundaries.
If two spatial surfaces with the same boundary are not equivalent,
they are not equivalent as Seifert surfaces.
Hence, 
distinguishing spatial surfaces with the same boundary play a role in addressing the classification problem of Seifert surfaces for the boundary.
A \textit{spatial trivalent graph} 
is a finite trivalent graph embedded in $S^{3}$.
As a remark,
we allow trivalent graphs to have loops, multiple edges, and circle components.
Diagrams of 
spatial trivalent graphs are defined in the same manner as in knot theory.
Every spatial surface can be presented by a diagram of a 
spatial trivalent graph, and
Reidemeister-type theorem for spatial surfaces was established \cite{Matsuzaki2021}.
This provides a tool for studying Seifert surfaces via their combinatorial structures.

\textit{Racks} \cite{Fenn-Rourke1992} and \textit{quandles} \cite{Joyce1982,Matveev1982} are algebraic structures 
whose axioms are corresponding to Reidemeister moves in knot theory.
A \textit{coloring} of a knot diagram by a given rack (resp. quandle) is a map 
from the set of all arcs of the diagram to the rack (resp. quandle)
satisfying a certain condition.
For each rack (resp. quandle),
the number of colorings of a diagram by the rack (resp. quandle) is an invariant of the framed knot (resp. knot) in $S^{3}$
presented by the diagram.
A \textit{multiple group rack} (MGR) \cite{Ishii-Matsuzaki-Murao2020} is a rack 
on a disjoint union of groups with a binary operation satifying some conditions 
which are motivated by Reidemeister moves for diagrams of spatial surfaces.
A \textit{coloring} of a diagram by an MGR is defined analogously to a coloring of a knot diagram by a rack.
For a given MGR,
the number of colorings of a diagram by the MGR is known to be an invariant of the spatial surface presented by the diagram.
  There are several examples of Seifert surfaces which are distinguished by the invariants \cite{Ishii-Matsuzaki-Murao2020}. 
  In this paper,
  we give an infinite family of pairs of Seifert surfaces for each link satisfying the following conditions:
  in each pair of Seifert surfaces, 
  (i) their regular neighborhoods are ambiently isotopic in $S^{3}$,
  (ii) their Seifert matrices are unimodularly congruent, and 
  (iii) they are not equivalent.

\section{Spatial surfaces and handlebody-knots}{\label{Sec:Spatial_surface}}

A \textit{Seifert surface} for a link is an 
oriented surface in the $3$-sphere $S^{3} = \mathbb{R}^{3} \sqcup \left\{\infty\right\}$ 
whose boundary is ambiently isotopic to the link.
It is well known that every link has a Seifert surface \cite{Frankl-Pontrjagin1930,Seifert1935}.
Two Seifert surfaces for a link are said to be \textit{equivalent} if they are ambiently isotopic in the exterior of the link.
A \textit{spatial surface} is a compact surface embedded in 
$S^{3}$. 
Two spatial surfaces $F_{1}$ and $F_{2}$ are said to be \textit{equivalent} (denoted by $F_{1} \cong F_{2}$) if they are ambiently isotopic in $S^{3}$.
When each component of a spatial surface is homeomorphic to an annulus,
then the surface is a framed link in $S^{3}$.
In this sense,
spatial surfaces can be regarded as a generalization of framed links in $S^{3}$.
Throughout this paper,
we assume that 
(1) a spatial surface is oriented and 
(2) each component of a spatial surface is neither a closed disk nor a surface without boundary.
Under the assumptions,
spatial surfaces are Seifert surfaces for their boundaries.
As a remark,
if two spatial surfaces with the same boundary are not equivalent,
then they are not equivalent as Seifert surfaces for the boundary. 
A \textit{spatial trivalent graph} is a finite trivalent graph embedded in $S^{3}$.
In this paper,
we allow trivalent graphs to have loops, multiple edges, and \textit{circle components}, i.e., edges without vertices.
Diagrams of spatial trivalent graphs are defined in the same manner as in knot theory.
An \textit{edge} of a diagram $D$ of a spatial trivalent graph refers to a sub-diagram of $D$ that presents an edge of the spatial trivalent graph.
In particular,
an edge of a diagram of a spatial trivalent graph that presents a circle component of the spatial trivalent graph is called a \textit{circle component} of the diagram.

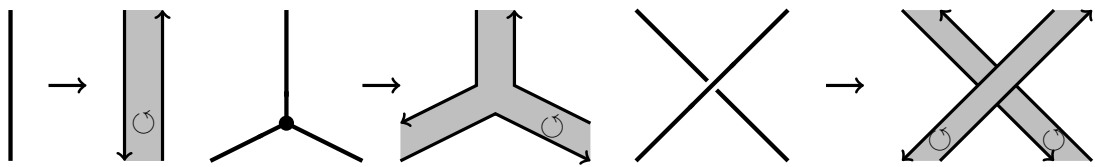
\begin{figure}[h]
  \centering
  \begin{tikzpicture}[line width=1.2pt, use Hobby shortcut]
    \begin{knot}[
      consider self intersections=true,
      ignore endpoint intersections=false,
      flip crossing/.list={}
    ]
    \strand(0,2)--(0,0);
    \end{knot}
    \fill[lightgray] (1.5,2)--(1.5,0)--(2,0)--(2,2);
    \draw[->](0.5,1)--(1,1);
    \draw[->](1.5,2)--(1.5,0);
    \draw[->](2,0)--(2,2);
    \draw (1.75,0.5) node{$\circlearrowleft $};
  \end{tikzpicture}
  \quad
  \begin{tikzpicture}[line width=1.2pt, use Hobby shortcut]
    \begin{knot}[
      consider self intersections=true,
      ignore endpoint intersections=false,
      flip crossing/.list={}
    ]
    \strand(1,2)--(1,0.5);
    \strand(1,0.5)--(0,0);
    \strand(1,0.5)--(2,0);
    \end{knot}
    \draw[->] (2,1)--(2.5,1);
    \fill[lightgray] (3.5,2)--(3.5,0.625)--(4,0.625)--(4,2);
    \fill[lightgray] (2.5,0)--(3.75,0.625)--(3.5,1)--(2.5,0.5);
    \fill[lightgray] (5,0)--(3.75,0.625)--(4,1)--(5,0.5);
    \fill (1,0.5) circle(0.1);
    \draw[->] (3.5,2)--(3.5,1)--(2.5,0.5);
    \draw[->] (2.5,0)--(3.75,0.625)--(5,0);
    \draw[->] (5,0.5)--(4,1)--(4,2);
    \draw (4.5,0.45) node{$\circlearrowleft$};
  \end{tikzpicture}
  \quad
  \begin{tikzpicture}[line width=1.2pt, use Hobby shortcut]
    \begin{knot}[
      consider self intersections=true,
      ignore endpoint intersections=false,
      flip crossing/.list={},
      only when rendering/.style={}
    ]
    \strand(2,2)--(0,0);
    \strand(0,2)--(2,0);
    \end{knot}
    \draw[->] (2.5,1)--(3,1);
    \fill[lightgray] (3.5,2)--(5.5,0)--(6,0)--(4,2);
    \draw[->] (3.5,2)--(5.5,0);
    \draw[->] (6,0)--(4,2);
    \fill[lightgray] (5.5,2)--(3.5,0)--(4,0)--(6,2);
    \draw[->] (5.5,2)--(3.5,0);
    \draw[->] (4,0)--(6,2);
    \draw (5.5,0.28) node{$\circlearrowleft$};
    \draw (4,0.28) node{$\circlearrowleft$};
  \end{tikzpicture}
  \caption{A process of constructing a spatial surface from a diagram}
  {\label{Fig:Construction_of_spatial_surfaces}}
\end{figure}

Let $D$ be a diagram of a spatial trivalent graph.
We construct from $D$ a spatial surface $F(D)$ as illustrated in Fig.~\ref{Fig:Construction_of_spatial_surfaces}.
Specifically,
consider a regular neighborhood $N(D) \subset \mathbb{R}^{2}$ of $D$ and 
replace each crossing with two bands embedded in $\mathbb{R}^{3}$ as depicted in the rightmost of Fig.~\ref{Fig:Construction_of_spatial_surfaces}.
Then the resulting surface is compact and embedded in $\mathbb{R}^{3}$, 
equipped with the orientation induced by that of $\mathbb{R}^{2}$.
Placing the oriented surface in $S^{3} = \mathbb{R}^{3} \sqcup \left\{\infty\right\},$ 
we obtain the spatial surface $F(D)$, which we call the \textit{spatial surface obtained from $D$}.
For any spatial surface $F$,
there exists a diagram $D$ of a spatial trivalent graph such that $F \cong F(D)$ \cite{Ishii-Matsuzaki-Murao2020,Matsuzaki2021}.
A \textit{diagram} of a spatial surface $F$ means a diagram $D$ of a spatial surface with $F(D) \cong F.$

In \cite{Matsuzaki2021},
a Reidemeister-type theorem for spatial surfaces was established.
We remark that a Reidemeister-type theorem for non-orientable spatial surfaces was established as well \cite{Matsuzaki2021}.
However, 
we do not consider the non-orientable case in the present paper.

\begin{figure}[h]
  \centering
  \begin{tikzpicture}[line width=1.6pt, use Hobby shortcut]
    \begin{knot}[
      consider self intersections=true,
      ignore endpoint intersections=false,
      flip crossing/.list={},
      only when rendering/.style={}
    ]
    \strand (1,2) -- (1,1.5) -- (0.5,0.5) -- (0.25,1) -- (0.5,1.5) -- (1,0.5) -- (1,0);
    \strand (3,2) -- (3,0);
    \strand (5.5,0) -- (5.5,0.5) -- (5,1.5) -- (4.75,1) -- (5,0.5) -- (5.5,1.5) -- (5.5,2);
    \end{knot}
    \draw[<->] (1.5,1)--(2.5,1);
    \draw (2,0.3) node{$\mathrm{R}1$};
    \draw[<->] (3.5,1)--(4.5,1);
    \draw (4,0.3) node{$\mathrm{R}1$};
   \end{tikzpicture}
  \ 
  \begin{tikzpicture}[line width=1.6pt, use Hobby shortcut]
    \begin{knot}[
      consider self intersections=true,
      ignore endpoint intersections=false,
      flip crossing/.list={},
      only when rendering/.style={}
    ]
    \strand (0,2)--(0,0);
    \strand (1,2)--(1,0);
    \strand (3,2)--(3.5,1.5)--(4,1)--(3.5,0.5)--(3,0);
    \strand (4,2)--(3.5,1.5)--(3,1)--(3.5,0.5)--(4,0);
    \end{knot}
    \draw[<->] (1.5,1)--(2.5,1);
    \draw (2,0.3) node{$\mathrm{R}2$};
   \end{tikzpicture}
  \ 
  \begin{tikzpicture}[line width=1.6pt, use Hobby shortcut]
    \begin{knot}[
      consider self intersections=true,
        ignore endpoint intersections=false,
        flip crossing/.list={},
        only when rendering/.style={}
      ]
      \strand(2,2)--(0,0);
      \strand(1,2)--(0.5,1.5)--(0,1)--(0.5,0.5)--(1,0);
      \strand(0,2)--(2,0);
      \strand(6,2)--(4,0);
      \strand(5,2)--(6,1)--(5,0);
      \strand(4,2)--(6,0);
      \end{knot}
      \draw[<->] (2.5,1)--(3.5,1);
      \draw (3,0.3) node{$\mathrm{R}3$};
    \end{tikzpicture}
    
    \vspace{5mm}
    \begin{tikzpicture}[line width=1.6pt, use Hobby shortcut]
      \begin{knot}[
        consider self intersections=true,
        ignore endpoint intersections=false,
        flip crossing/.list={},
        only when rendering/.style={}
      ]
      \strand (2.25,1.5) -- (2.25,1.125) -- (1.75,0.375) -- (1.5,0.75) -- (1.75,1.125) -- (2.25,0.375) -- (2.25,0);
      \strand (10.75,0) -- (10.75,0.375) -- (10.25,1.125) -- (10,0.75) -- (10.25,0.375) -- (10.75,1.125) -- (10.75,1.5);
      \end{knot}
      \draw (0,0.75) -- (1.5,0.75);
      \filldraw (1.5,0.75) circle(0.1);
      \draw[<->] (2.75,0.75)--(3.75,0.75);
      \draw (3.25,0.225) node{$\mathrm{R}4$};
      \draw (4.25,0.75) -- (5.75,0.75);
      \draw (6.5,1.5) -- (5.75,0.75) -- (6.5,0);
      \filldraw (5.75,0.75) circle(0.1);
      \draw[<->] (7,0.75) -- (8,0.75);
      \draw (7.5,0.225) node{$\mathrm{R}4$};
      \draw (8.5,0.75) -- (10,0.75);
      \filldraw (10,0.75) circle(0.1);
    \end{tikzpicture}

    \vspace{5mm}
    \begin{tikzpicture}[line width=1.6pt, use Hobby shortcut]
      \begin{knot}[
        consider self intersections=true,
        ignore endpoint intersections=false,
        flip crossing/.list={},
        only when rendering/.style={}
      ]
      \strand(0.75,1.5)--(0.75,0);
      \strand(0,1.5)--(1.5,0.75)--(0,0);
      \strand(1.5,0.75)--(3,0.75);
      \strand(4.5,1.5)--(6,0.75)--(4.5,0);
      \strand(6.75,1.5)--(6.75,0);
      \strand(6,0.75)--(7.5,0.75);
      \end{knot}
      \filldraw (1.5,0.75) circle(0.1);
      \draw[<->] (3.375,0.75)--(4.125,0.75);
      \filldraw (6,0.75) circle(0.1);
      \draw (3.75,0.225) node{R$5$};
    \end{tikzpicture}
		\quad
		\begin{tikzpicture}[line width=1.6pt, use Hobby shortcut]
      \begin{knot}[
        consider self intersections=true,
        ignore endpoint intersections=false,
        flip crossing/.list={},
        only when rendering/.style={}
      ]
      \strand(0,1.5)--(1.5,0.75)--(0,0);
      \strand(0.75,1.5)--(0.75,0);
      \strand(1.5,0.75)--(3,0.75);
      \strand(4.5,1.5)--(6,0.75)--(4.5,0);
      \strand(6,0.75)--(7.5,0.75);
      \strand(6.75,1.5)--(6.75,0);
      \end{knot}
      \filldraw (1.5,0.75) circle(0.1);
      \draw[<->] (3.375,0.75)--(4.125,0.75);
      \filldraw (6,0.75) circle(0.1);
      \draw (3.75,0.225) node{R$5$};
    \end{tikzpicture}

        \vspace{5mm}
    \begin{tikzpicture}[line width=1.6pt, use Hobby shortcut]
      \begin{knot}[
        consider self intersections=true,
        ignore endpoint intersections=false,
        flip crossing/.list={},
        only when rendering/.style={}
      ]
      \strand(0,2)--(1,1.5)--(2,2);
      \strand(1,1.5)--(1,0.5);
      \strand(0,0)--(1,0.5)--(2,0);
      \strand(6,2)--(6.5,1)--(6,0);
      \strand(6.5,1)--(7.5,1);
      \strand(8,0)--(7.5,1)--(8,2);
      \end{knot}
      \filldraw (1,1.5) circle(0.1);
      \filldraw (1,0.5) circle(0.1);
      \draw[<->] (3,1)--(5,1);
      \draw (4,0.3) node{R$6$};
      \filldraw (6.5,1) circle(0.1);
      \filldraw (7.5,1) circle(0.1);
    \end{tikzpicture}
    \caption{Local moves on diagrams of spatial trivalent graphs}{\label{Fig:R-moves}}
\end{figure}
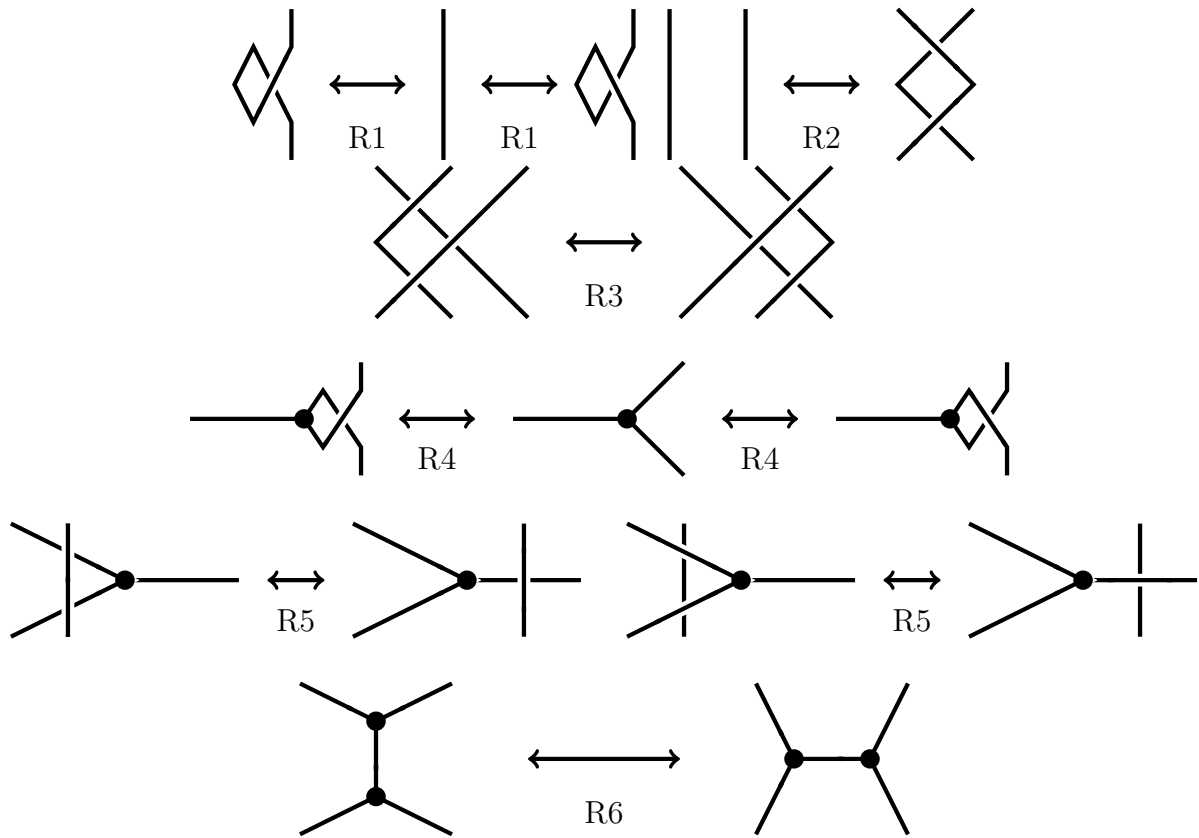

\thm[\cite{Matsuzaki2021}]{\label{Thm:R-moves}}{
  Two spatial surfaces are equivalent if and only if their diagrams are related by 
  a finite sequence of $\mathrm{R}2, \mathrm{R}3, \mathrm{R}5, \mathrm{R}6$ moves, as illustrated in Fig.~\ref{Fig:R-moves},
  and isotopies in $S^{2}$.
}\upshape

The set of the four moves 
$\mathrm{R}2$, $\mathrm{R}3$, $\mathrm{R}5$, and $\mathrm{R}6$ is referred to as 
the \textit{Reidemeister moves for spatial surfaces}.




A \textit{handlebody-knot} \cite{Ishii2008} is a handlebody embedded in $S^{3}$.
Two handlebody-knots $H_{1}$ and $H_{2}$ are said to be \textit{equivalent} ($H_{1} \cong H_{2}$) if
they are ambiently isotopic in $S^{3}$.
Handlebody-knots can be obtained as regular neighborhoods of spatial trivalent graphs in $S^{3}$.
A \textit{diagram} of a handlebody-knot $H$ is defined as 
a diagram of a spatial trivalent graph whose regular neighborhood is equivalent to $H$.
Let $D$ be a diagram of a spatial trivalent graph $G$.
We denote by $H(D)$ the handlebody-knot whose diagram is $D$.
In \cite{Ishii2008}, a Reidemeister-type theorem for handlebody-knots was introduced.

\thm[\cite{Ishii2008}]{\label{Thm:R-moves_for_Hdlebody-knots}}{
  Two handlebody-knots are equivalent if and only if their diagrams are related by a finite sequence of 
  $\mathrm{R}1$--$\mathrm{R}6$ moves, depicted in Fig.~\ref{Fig:R-moves}, and isoopies in $S^{2}$.
}\upshape

We refer to the set of moves $\mathrm{R}1$--$\mathrm{R}6$ as the \textit{Reidemeister moves for handlebody-knots}.

\rem{
  The Reidemeister moves for spatial surfaces are obtained from those for handlebody-knots 
  by removing the moves $\mathrm{R}1$ and $\mathrm{R}4$. 
  For diagrams of spatial surfaces, 
  the moves $\mathrm{R}1$ and $\mathrm{R}4$ change the framing of the spatial surfaces. 
  Therefore,
  applying the moves $\mathrm{R}1$ and $\mathrm{R}4$ is not allowed for diagrams of spatial surfaces.
}

\section{Seifert matrices}{\label{Sec:Seifert_matrix}}

We denote the linking number of disjoint oriented knots $K_{1}$ and $K_{2}$ by $\mathrm{Lk}(K_{1}, K_{2})$.
  Let $F$ be a spatial surface.
  For a simple closed curve $l$ on $F$,
  we denote by $l^{+}$ a parallel copy of $l$
  obtained by pushing $l$ in the direction of the positive normal orientation of $F$.
  Take the bilinear map $\phi: H_{1}(F) \times H_{1}(F) \to \mathbb{Z}$, defined by $\phi(x, y) = \mathrm{Lk}(l_{x}^{+},l_{y})$
  where $l_{x}$ and $l_{y}$ are simple closed curves on $F$ representing the homology classes $x$ and $y$.
  Given a basis of $H_{1}(F)$,
  a presentation matrix $V = V(F)$ of $\phi$ with respect to the basis is called a \textit{Seifert matrix} of $F$.
A \textit{unimodular matrix} is a square matrix over $\mathbb{Z}$ with determinant $\pm 1$.
Two square matrices $V_{1}$ and $V_{2}$ over $\mathbb{Z}$ are said to be \textit{unimodularly congruent} 
if there exists a unimodular matrix $P$ such that $V_{2} = P^{T} V_{1} P$,
where $P^{T}$ is the transpose of $P$.

Although a Seifert matrix of a spatial surface depends on the choice of a basis for the first homology group of the spatial surface,
the following result holds for Seifert matrices for spatial surfaces.

\prop{\label{Prop:Seifert_matrix}}{
  If two spatial surfaces $F_{1}$ and $F_{2}$ are equivalent,
  then their Seifert matrices $V(F_{1})$ and $V(F_{2})$ are unimodularly congruent.
}\upshape



According to Proposition~\ref{Prop:Seifert_matrix},
a Seifert matrix of a spatial surface is an invariant of the spatial surface up to unimodularly congruent.



\cor{\label{Cor:Seifert_matrix}}{
  Let $V_{1}$ and $V_{2}$ be Seifert matrices of spatial surfaces $F_{1}$ and $F_{2}$, respectively.
  If $F_{1}$ and $F_{2}$ are equivalent,
  then $\gcd \left\{k \times k\mbox{-minors of}\ V_{1}\right\} = \gcd \left\{k \times k\mbox{-minors of}\ V_{2}\right\}$ 
  for any $k \in \left\{1, 2, \ldots, l\right\}$, 
  where $l$ is the size of Seifert matrices $V_{1}$ and $V_{2}$.
}\upshape

\section{Multiple group racks}{\label{Sec:MGRs}}

\defn[\cite{Fenn-Rourke1992}]{\label{Def:Rack}}{
  Let $R$ be a nonempty set and $\triangleleft : R \times R \to R$ be a binary operation on $R$.
  The pair $(R, \triangleleft)$ is a \textit{rack} if it satisfies the following conditions (i)--(ii).
  \begin{itemize}
    \item[(i)] For any $y \in R$, the map $S_{y}: R \to R$, defined by $S_{y}(x) = x \triangleleft y$, is bijective.
    \item[(ii)] For any $x, y, z \in R$, $(x \triangleleft y) \triangleleft z = (x \triangleleft z) \triangleleft (y \triangleleft z)$.
  \end{itemize}
  A \textit{quandle} \cite{Joyce1982,Matveev1982} is a rack $(R, \triangleleft)$ such that it satisfies $x \triangleleft x = x$ for any $x \in R$.
}

\example{\label{Ex:Racks_and_quandles}}{
  \begin{itemize}
    \item Let $n$ be a positive integer and $\triangleleft: \mathbb{Z}_{n} \times \mathbb{Z}_{n} \to \mathbb{Z}_{n}$ the binary operation on $\mathbb{Z}_{n}$, defined by $x \triangleleft y = 2y - x$.
    Then, the pair $R_{n} = (\mathbb{Z}_{n}, \triangleleft)$ is a quandle.
    The quandle $R_{n}$ is called the \textit{dihedral quandle}.
    \item Let $n$ be a positive integer and $\triangleleft: \mathbb{Z}_{n} \times \mathbb{Z}_{n} \to \mathbb{Z}_{n}$ the binary operation on $\mathbb{Z}_{n}$, defined by $x \triangleleft y = x + 1$.
    Then, the pair $C_{n} = (\mathbb{Z}_{n}, \triangleleft)$ is a rack.
    The rack $R_{n}$ is called the \textit{cyclic rack}.
    \item Let $(R, \triangleleft), (R', \triangleleft')$ be racks and 
    $\ast: (R \times R') \times (R \times R') \to R \times R'$ the binary operation on $R \times R'$, defined by $(x, x') \ast (y, y') = (x \triangleleft y, x' \triangleleft' y')$.
    Then,
    the pair $(R \times R', \ast)$ is a rack.
    The rack $(R \times R', \ast)$ is called the \textit{product} of the racks $(R, \triangleleft)$ and $(R', \triangleleft')$.
  \end{itemize}
}

\defn[\cite{Ishii-Matsuzaki-Murao2020}]{\label{Def:MGR}}{
  Let $X = \bigsqcup_{\lambda \in \Lambda} G_{\lambda}$ be a disjoint union of groups $G_{\lambda}$ ($\lambda \in \Lambda$) with the identity elements $e_{\lambda}$ and 
  $\triangleleft : X \times X \to X$ a binary operation on $X$.
  The pair $(X, \triangleleft)$ is a \textit{multiple group rack} (MGR) if it satisfies the following conditions (i)--(iii).
  \begin{itemize}
    \item[(i)] For any $x \in X$, $\lambda \in \Lambda$, and $a, b \in G_{\lambda}$, $x \triangleleft (ab) = (x \triangleleft a) \triangleleft b$ and $x \triangleleft e_{\lambda} = x$.
    \item[(ii)] For any $x, y, z \in X$, $(x \triangleleft y) \triangleleft z = (x \triangleleft z) \triangleleft (y \triangleleft z)$.
    \item[(iii)] For any $x \in X$ and $\lambda \in \Lambda$, there exists $\mu \in \Lambda$ such that
    for any $a, b \in X$, $a \triangleleft x, b \triangleleft x \in G_{\mu}$ and $(ab) \triangleleft x = (a \triangleleft x) (b \triangleleft x)$.
  \end{itemize}
  A \textit{multiple conjugation quandle} (MCQ) \cite{Ishii2015} is an MGR $(X = \bigsqcup_{\lambda \in \Lambda} G_{\lambda}, \triangleleft)$ such that
  for any $\lambda \in \Lambda$ and $a, b \in G_{\lambda}$, $a \triangleleft b = b^{-1} a b$.
}

When there is no risk of confusion,
we sometimes refer an MGR $(X, \triangleleft)$ (resp. a rack $(R, \triangleleft)$) simply as $X$ (resp. $R$).

\prop{\label{Prop:MGR_is_rack}}{
  Let $(X, \triangleleft)$ be an MGR.
  Then, $(X, \triangleleft)$ is a rack.
}\upshape

\begin{proof}
  We show that for any $y \in X$, 
  the map $S_{y}: X \times X \to X$, defined by $S_{y}(x) = x \triangleleft y$, 
  is bijective.
  Since $X$ is a disjoint union of groups,
  there is a unique group $G$ such that $y \in G$.
  Take the inverse element $y^{-1}$ of $y$ with respect to the group operation of $G$.
  By the condition (i) of Definition~\ref{Def:MGR},
  \begin{equation*}
    (x \triangleleft y) \triangleleft y^{-1} = x \triangleleft (y y^{-1}) = x \triangleleft e = x,
  \end{equation*}
  where $e$ is the identity element of the group $G$.
  Thus,
  the map $S_{y^{-1}}: X \times X \to X$, defined by $S_{y^{-1}}(x) = x \triangleleft y^{-1}$, 
  is the inverse map of $S_{y}$.
  Hence, 
  the map $S_{y}$ is bijective.

  Since the condition (ii) of Definition~\ref{Def:Rack} coincides the condition (ii) of Definition~\ref{Def:MGR},
  we conclude that the MGR $(X, \triangleleft)$ is a rack.
\end{proof}

\example{\label{Ex:MGRs}}{
  Let $R$ be a finite rack.
  Define $n = \min \left\{k \in \mathbb{Z}_{>0} \mid \mbox{for any}\ x, y \in R, S_{y}^{k}(x) = x\right\}$.
  Then, 
  $R \times \mathbb{Z}_{n} = \bigsqcup_{x \in R} \left(\left\{x\right\} \times \mathbb{Z}_{n}\right)$ is an MGR with the following operations.
  \begin{gather*}
    (x, i) \triangleleft (y, j) = (S_{y}^{j}(x), i), \quad (x, i)(x, j) = (x, i+j).
  \end{gather*}
  We observe that this MGR $R \times \mathbb{Z}_{n}$ is the disjoint union of copies of the group $\mathbb{Z}_{n}$, 
  indexed by the elements of $R$.
}

\begin{figure}[h]
  \centering
  \begin{tikzpicture}[line width=2pt]
    \draw (0,2)--(1,1)--(1,0);
    \draw (2,2)--(1,1);
    \draw (0.5,1.5) node{\rotatebox{315}{$\blacktriangleright$}};
    \draw (1.5,1.5) node{\rotatebox{225}{$\blacktriangleright$}};
    \draw (1,0.5) node{\rotatebox{270}{$\blacktriangleright$}};
    \filldraw (1,1) circle(0.125);
    \draw (1,-0.5) node{Y-orientation};
  \end{tikzpicture}
  \quad 
  \begin{tikzpicture}[line width=2pt]
    \draw (0,2)--(1,1)--(1,0);
    \draw (2,2)--(1,1);
    \draw (0.5,1.5) node{\rotatebox{135}{$\blacktriangleright$}};
    \draw (1.5,1.5) node{\rotatebox{45}{$\blacktriangleright$}};
    \draw (1,0.5) node{\rotatebox{90}{$\blacktriangleright$}};
    \filldraw (1,1) circle(0.125);
    \draw (1,-0.5) node{Y-orientation};
  \end{tikzpicture}
  \quad
  \begin{tikzpicture}[line width=2pt]
    \draw (0,2)--(1,1)--(1,0);
    \draw (2,2)--(1,1);
    \draw (0.5,1.5) node{\rotatebox{315}{$\blacktriangleright$}};
    \draw (1.5,1.5) node{\rotatebox{235}{$\blacktriangleright$}};
    \draw (1,0.5) node{\rotatebox{90}{$\blacktriangleright$}};
    \filldraw (1,1) circle(0.125);
    \draw (1,-0.5) node{sink};
  \end{tikzpicture}
  \quad
  \begin{tikzpicture}[line width=2pt]
    \draw (0,2)--(1,1)--(1,0);
    \draw (2,2)--(1,1);
    \draw (0.5,1.5) node{\rotatebox{135}{$\blacktriangleright$}};
    \draw (1.5,1.5) node{\rotatebox{45}{$\blacktriangleright$}};
    \draw (1,0.5) node{\rotatebox{270}{$\blacktriangleright$}};
    \filldraw (1,1) circle(0.125);
    \draw (1,-0.5) node{source};
  \end{tikzpicture}
  \caption{All orientations around trivalent vertices}{\label{Fig:All_orientation}}
\end{figure}
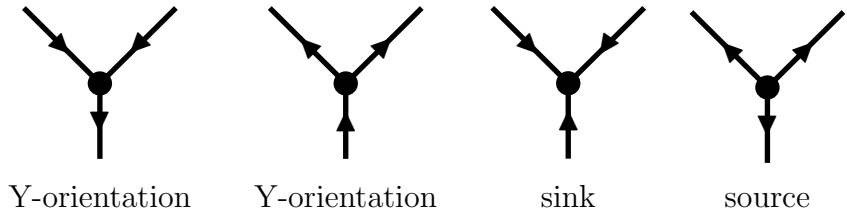

Let $D$ be a diagram of a spatial surface.
A \textit{Y-orientation} of $D$ is an assignment of orientations to all edges 
such that there are no sources and sinks as illustrated in Fig.~\ref{Fig:All_orientation}.
Every diagram of an spatial trivalent graph admits at least one Y-orientation \cite{Ishii2015,Lebed2015}.
A diagram with a Y-orientation is called a \textit{Y-oriented diagram}.
We remark that there is no relation between a Y-orientation of a diagram $D$ and the orientation of the spatial surface $F(D)$.
An \textit{arc} of $D$ is a connected component of $D$ obtained by cutting $D$ at undercrossings and vertices.
We denote by $\mathcal{A}(D)$ the set of all arcs of $D$.

\defn{\label{Def:MGRCols}}{
Let $X = \bigsqcup_{\lambda \in \Lambda} G_{\lambda}$ be an MGR and
$D$ a Y-oriented diagram of a spatial surface.
An \textit{$X$-coloring} of a Y-oriented diagram $D$ or a \textit{coloring} of $D$ by $X$ is a map $C: \mathcal{A}(D) \to X$ satisfying the following conditions.
\begin{itemize}
  \item For each crossing of $D$, $C$ satisfies $C(a_{i}) \triangleleft C(a_{j}) = C(a_{k})$,
  where $a_{i}, a_{j}, a_{k} \in \mathcal{A}(D)$ are as shown in the left side of Fig.~\ref{Fig:X-coloring_conditions}.
  \item For each vertex of $D$, $C$ satisfies the condition that $C(a_{i}), C(a_{j})$ and $C(a_{k})$ are in the same group, and that
  $C(a_{i}) C(a_{j}) = C(a_{k})$ holds, where $a_{i}, a_{j}, a_{k} \in \mathcal{A}(D)$ are as shown in the center or the right side of Fig.~\ref{Fig:X-coloring_conditions}.  
\end{itemize}
\begin{figure}[h]
  \centering
    \begin{tikzpicture}[line width=2pt, use Hobby shortcut]
      \begin{knot}[
        consider self intersections=true,
        ignore endpoint intersections=false,
        flip crossing/.list={}
        ]
      \strand (2,2)--(0,0);
      \strand (0,2)--(2,0);
      \end{knot}  
      \draw (0,0) node{\rotatebox{225}{$\blacktriangleright$}};
      \draw (-0.3,2) node{$a_{i}$};
      \draw (2.3,2) node{$a_{j}$};
      \draw (2.3,0) node{$a_{k}$};
      \draw (1,-0.5) node{$C(a_{i}) \triangleleft C(a_{j}) = C(a_{k})$};
    \end{tikzpicture}
    \quad
    \begin{tikzpicture}[line width=2pt]
      \draw (0,2)--(1,1)--(1,0);
      \draw (2,2)--(1,1);
      \draw (0.5,1.5) node{\rotatebox{315}{$\blacktriangleright$}};
      \draw (1.5,1.5) node{\rotatebox{225}{$\blacktriangleright$}};
      \draw (1,0.5) node{\rotatebox{270}{$\blacktriangleright$}};
      \draw (-0.3,2) node{$a_{i}$};
      \draw (2.3,2) node{$a_{j}$};
      \draw (1.5,0) node{$a_{k}$};
      \filldraw (1,1) circle(0.125);
      \draw (1,-0.5) node{$C(a_{i})C(a_{j}) = C(a_{k})$};
    \end{tikzpicture}
    \quad 
    \begin{tikzpicture}[line width=2pt]
      \draw (1,2)--(1,1)--(0,0);
      \draw (2,0)--(1,1);
      \draw (0.5,0.5) node{\rotatebox{225}{$\blacktriangleright$}};
      \draw (1.5,0.5) node{\rotatebox{315}{$\blacktriangleright$}};
      \draw (1,1.5) node{\rotatebox{270}{$\blacktriangleright$}};
      \filldraw (1,1) circle(0.125);
      \draw (1.5,2) node{$a_{k}$};
      \draw (-0.3,0) node{$a_{i}$};
      \draw (2.3,0) node{$a_{i}$};
      \draw (1,-0.5) node{$C(a_{i})C(a_{j}) = C(a_{k})$};
    \end{tikzpicture} 
    \caption{$X$-coloring conditions}{\label{Fig:X-coloring_conditions}}  
\end{figure}
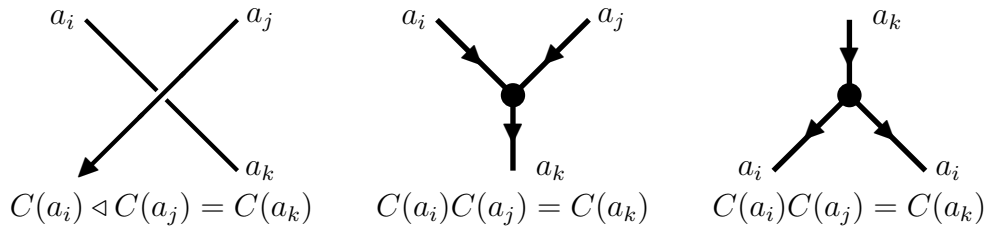
We denote the set of all $X$-colorings of $D$ by $\mathrm{Col}_{X}(D)$.
}

\thm[\cite{Ishii-Matsuzaki-Murao2020}]{\label{Thm:Col_number}}{
  Let $X$ be an MGR and
  $D$ and $D'$ Y-oriented diagrams of spatial surfaces.
  If $F(D) \cong F(D')$,
  then there exists a bijection between $\mathrm{Col}_{X}(D)$ and $\mathrm{Col}_{X}(D')$.
  In particular,
  the cardinality $\left\lvert \mathrm{Col}_{X}(D)\right\rvert$ is an invariant of the spatial surface $F(D)$.
}\upshape

Note that the cardinality $\left\lvert \mathrm{Col}_{X}(D)\right\rvert$ in Theorem~\ref{Thm:Col_number} is 
independent of the choice of a Y-orientation of $D$.

\section{Infinitely many pairs of spatial surfaces}{\label{Sec:Main}}

\thm{\label{Thm:Main_result}}{
  For any oriented link $L$, 
  there is an infinitely many pairs $\left\{(F_{n}, {F_{n}}')\right\}_{n \in \mathbb{Z}}$ of Seifert surfaces for $L$ such that 
  the following statements (i)--(iii).
  \begin{itemize}
    \item[(i)] For any $n \in \mathbb{Z}$, the regular neighborhoods $N(F_{n})$ and $N({F_{n}}')$ of $F_{n}$ and ${F_{n}}'$ are equivalent as handlebody-knots.
    \item[(ii)] For any $n \in \mathbb{Z}$, Seifert matrices of $F_{n}$ and ${F_{n}}'$ are unimodularly congruent.
    \item[(iii)] For any $n \in \mathbb{Z}$, $F_{n} \not\cong {F_{n}}'$.
  \end{itemize}
}\upshape



\begin{proof}
  (1) We give a family $\left\{(F_{n}, {F_{n}}')\right\}_{n \in \mathbb{Z}}$ of pairs of Seifert surfaces for $L$, as follows.
  
  Take a diagram $D$ of a Seifert surface $F$ for $L$
  as illustrated in Fig.~\ref{Fig:Diagram}.
  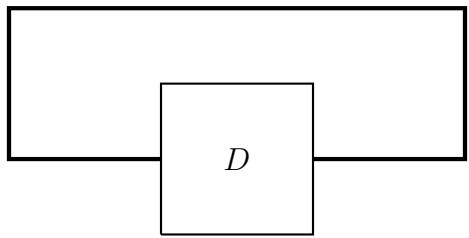
\begin{figure}[h]
    \centering
    \begin{tikzpicture}[line width=1.6pt]
      \draw[line width=0.8pt] (2,0) -- (4,0) -- (4,2) -- (2,2) -- (2,0);
      \draw (3,1) node{$D$};
      \draw (2,1) -- (0,1) -- (0,3) -- (6,3) -- (6,1) -- (4,1);
    \end{tikzpicture} 
    \caption{A diagram $D$ of the Seifert surface $F$ for $L$}{\label{Fig:Diagram}}
  \end{figure}
  For the diagram $D$, 
  attach edges in two different ways as illustrated in Fig.~\ref{Fig:Diagrams}, 
  and denote the resulting diagrams by $D_{n}$ and $D_{n}'$, where $n \in \mathbb{Z}$ is an integer.
  \begin{figure}[h]
    \centering
    \begin{tikzpicture}[line width=1.6pt, scale=0.8]
      \draw (7,7) -- (7,2);
      \draw[white, line width=3.8pt] (6,6) -- (8,6);\
      \draw (7,7) -- (8,7) -- (8,3) -- (7,3);

      \draw (3,-1) node{$D_{n}$};
      \draw (4.5,6) -- (8,6);
      \draw (3.5,6) -- (2.5,5);
      \draw[white, line width=3.8pt] (3.25,5.25) -- (2.75,5.75);
      \draw (4.5,6) -- (3.5,5);
      \draw[white, line width=3.8pt] (4.25,5.25) -- (3.75,5.75);
      \draw (2.5,6) -- (1.5,5);
      \draw[white, line width=3.8pt] (2.25,5.25) -- (1.75,5.75);
      \draw (4.5,5) -- (3.5,6);
      \draw (3.5,5) -- (2.5,6);
      \draw (2.5,5) -- (1.5,6);
      \draw (1.5,5) -- (1.5,4) -- (4.5,4) -- (4.5,5);
      \draw (1.5,6) -- (1.25,6);
      \draw (1.25,6.5) -- (0.75,6);
      \draw[white, line width=3.8pt] (0.95,6.2) -- (1.05,6.3);
      \draw (1.25,6) -- (0.75,6.5) -- (0.75,6.75) -- (1.25,6.75) -- (1.25,6.5);
      \draw (0.75,6) -- (-1,6);

      \draw (4.25,7) -- (7,7);
      \draw (4.25,7.5) -- (3.75,7);
      \draw[white, line width=3.8pt] (3.95,7.2) -- (4.05,7.3);
      \draw (4.25,7) -- (3.75,7.5) -- (3.75,7.75) -- (4.25,7.75) -- (4.25,7.5);
      \draw (3.75,7) -- (3.25,7);
      \draw (3.25,7.5) -- (2.75,7);
      \draw[white, line width=3.8pt] (2.95,7.2) -- (3.05,7.3);
      \draw (3.25,7) -- (2.75,7.5) -- (2.75,7.75) -- (3.25,7.75) -- (3.25,7.5);
      \draw (2.75,7) -- (2.25,7);
      \draw (2.25,7.5) -- (1.75,7);
      \draw[white, line width=3.8pt] (1.95,7.2) -- (2.05,7.3);
      \draw (2.25,7) -- (1.75,7.5) -- (1.75,7.75) -- (2.25,7.75) -- (2.25,7.5);
      \draw (1.75,7) -- (1.25,7);
      \draw (1.25,7.5) -- (0.75,7);
      \draw[white, line width=3.8pt] (0.95,7.2) -- (1.05,7.3);
      \draw (1.25,7) -- (0.75,7.5) -- (0.75,7.75) -- (1.25,7.75) -- (1.25,7.5);
      \draw[line width=0.8pt] (-0.125,7) circle(0.5);
      \draw (-0.125,7) node{$2n$};
      \draw (0.75,7) -- (0.375,7);
      \draw (-0.625,7) -- (-1,7);

      \draw (-1,7) -- (-1,2) -- (0,2);
      
      \draw (5,2) -- (7,2);
      \draw[white,line width=3.8pt] (6,3) -- (6,1);
      \draw[line width=0.8pt] (2,0) -- (4,0) -- (4,2) -- (2,2) -- (2,0);
      \draw (3,1) node{$D$};
      \draw (2,1) -- (0,1) -- (0,3) -- (6,3) -- (6,1) -- (4,1);
      \draw (5,3) -- (5,1);
    \end{tikzpicture} 
    \ 
    \begin{tikzpicture}[line width=1.6pt, scale=0.8]
      \draw (7,7) -- (7,2);
      \draw[white, line width=3.8pt] (6,4) -- (8,4);
      \draw (7,7) -- (8,7) -- (8,3) -- (7,3);

      \draw (3,-1) node{$D_{n}'$};
      \draw (4.5,7) -- (8,7);
      \draw (3.5,7) -- (2.5,6);
      \draw[white, line width=3.8pt] (3.25,6.25) -- (2.75,6.75);
      \draw (4.5,7) -- (3.5,6);
      \draw[white, line width=3.8pt] (4.25,6.25) -- (3.75,6.75);
      \draw (2.5,7) -- (1.5,6);
      \draw[white, line width=3.8pt] (2.25,6.25) -- (1.75,6.75);
      \draw (4.5,6) -- (3.5,7);
      \draw (3.5,6) -- (2.5,7);
      \draw (2.5,6) -- (1.5,7);
      \draw (1.5,6) -- (1.5,5) -- (4.5,5) -- (4.5,6);
      \draw (1.5,7) -- (1.25,7);
      \draw (1.25,7.5) -- (0.75,7);
      \draw[white, line width=3.8pt] (0.95,7.2) -- (1.05,7.3);
      \draw (1.25,7) -- (0.75,7.5) -- (0.75,7.75) -- (1.25,7.75) -- (1.25,7.5);
      \draw[line width=0.8pt] (-0.125,7) circle(0.5);
      \draw (-0.125,7) node{$2n$};
      \draw (0.75,7) -- (0.375,7);
      \draw (-0.625,7) -- (-1,7);

      \draw (4.25,4) -- (8,4);
      \draw (4.25,4.5) -- (3.75,4);
      \draw[white, line width=3.8pt] (3.95,4.2) -- (4.05,4.3);
      \draw (4.25,4) -- (3.75,4.5) -- (3.75,4.75) -- (4.25,4.75) -- (4.25,4.5);
      \draw (3.75,4) -- (3.25,4);
      \draw (3.25,4.5) -- (2.75,4);
      \draw[white, line width=3.8pt] (2.95,4.2) -- (3.05,4.3);
      \draw (3.25,4) -- (2.75,4.5) -- (2.75,4.75) -- (3.25,4.75) -- (3.25,4.5);
      \draw (2.75,4) -- (2.25,4);
      \draw (2.25,4.5) -- (1.75,4);
      \draw[white, line width=3.8pt] (1.95,4.2) -- (2.05,4.3);
      \draw (2.25,4) -- (1.75,4.5) -- (1.75,4.75) -- (2.25,4.75) -- (2.25,4.5);
      \draw (1.75,4) -- (1.25,4);
      \draw (1.25,4.5) -- (0.75,4);
      \draw[white, line width=3.8pt] (0.95,4.2) -- (1.05,4.3);
      \draw (1.25,4) -- (0.75,4.5) -- (0.75,4.75) -- (1.25,4.75) -- (1.25,4.5);
      \draw (0.75,4) -- (-1,4);
      \draw (-1,7) -- (-1,2) -- (0,2);
      
      \draw (5,2) -- (7,2);
      \draw[white,line width=3.8pt] (6,3) -- (6,1);
      \draw[line width=0.8pt] (2,0) -- (4,0) -- (4,2) -- (2,2) -- (2,0);
      \draw (3,1) node{$D$};
      \draw (2,1) -- (0,1) -- (0,3) -- (6,3) -- (6,1) -- (4,1);
      \draw (5,3) -- (5,1);
    \end{tikzpicture} 
    
    \begin{tikzpicture}[line width=1.6pt]
      \draw (-6,4) -- (-6,3.5);
      \draw (-6,2.5) -- (-6,2);
      \draw[line width=0.8pt] (-6,3) circle(0.5);
      \draw (-6,3) node{$n$};
      \draw (0,3) node{: positive $n$ kinks (if $n \geq 0$) or negative $\left\lvert n\right\rvert$ kinks (if $n < 0$)};
    \end{tikzpicture}
    \caption{Diagrams $D_{n}$ and $D_{n}'$ ($n \in \mathbb{Z}$)}{\label{Fig:Diagrams}}
  \end{figure}
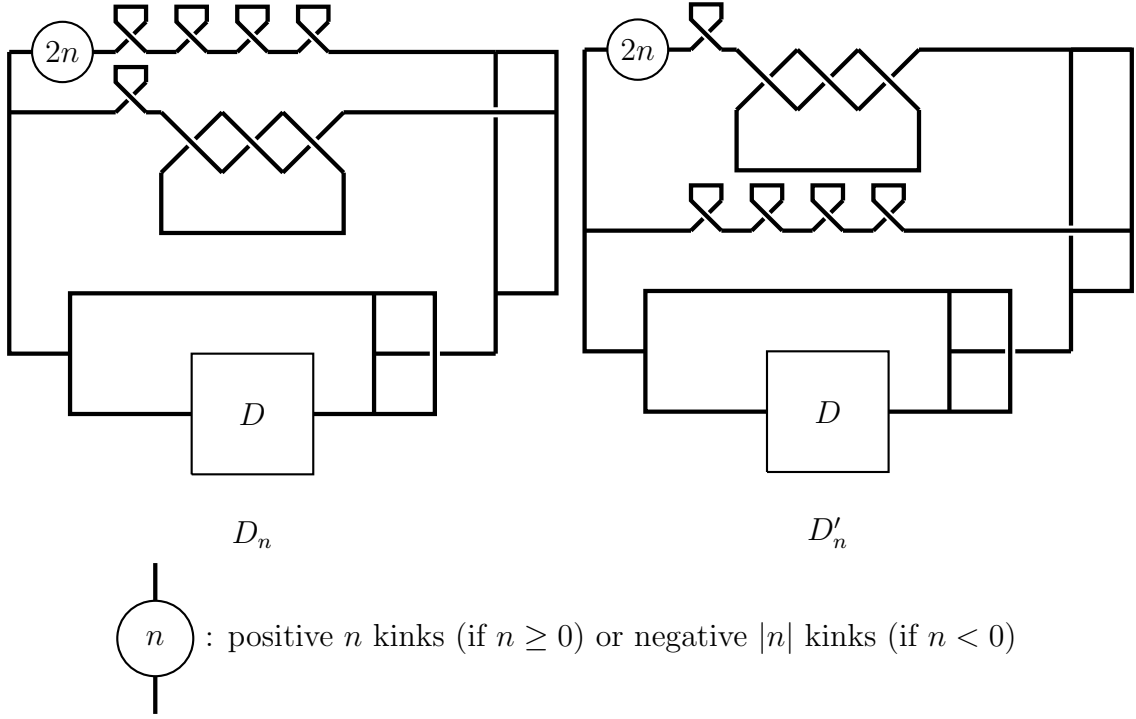
  For each $n \in \mathbb{Z}$, 
  we define $F_{n}$ and $F_{n}'$ to be the spatial surfaces $F(D_{n})$ and $F(D_{n}')$, respectively. 
  The links $\partial F(D_{n})$ and $\partial F(D_{n}')$ are ambiently isotopic.
  Hence, for any $n \in \mathbb{Z}$,
  these surfaces $F_{n}$ and $F_{n}'$ are Seifert surfaces for $L$.
  Thus,
  we obtain the desired family $\left\{(F_{n}, {F_{n}}')\right\}_{n \in \mathbb{Z}}$.
  
  (2) We show that the family $\left\{(F_{n}, {F_{n}}')\right\}_{n \in \mathbb{Z}}$ contains infinitely many mutually distinct pairs.



  Let $V$ be an $l \times l$ Seifert matrix of the spatial surface $F(D)$.
  Then,
  for each $k \in \mathbb{Z}$,
  $F_{k}$ has the Seifert matrix $V_{k} = 
  \begin{pmatrix}
    0 & 0 & 0 & 0 & \mathbf{0} \\
    0 & 0 & 1 & 0 & \mathbf{0} \\
    -1 & 0 & 4 & 0 & \mathbf{0} \\
    -1 & 0 & 0 & 4 + 2k & \mathbf{0} \\
    \mathbf{0} & \mathbf{0} & \mathbf{0} & \mathbf{0} & V
  \end{pmatrix}$.
  By Corollary~\ref{Cor:Seifert_matrix},
  for any $m, n \in \mathbb{Z}$ with $\left\lvert 4 + 2m\right\rvert \neq \left\lvert 4 + 2n\right\rvert$,
  it is enough to verify that 
  the greatest common divisor $\gcd \left\{k \times k\mbox{-minors of}\ V_{m}\right\}$ 
  is not equal to $\left\{l \times l\mbox{-minors of}\ V_{n}\right\}$ for some $l \in \mathbb{Z}$.
  Consider the integer 
  $s = \max \left\{i \in \mathbb{Z} \mid \mbox{there exists a non-zero}\  (i \times i)\mbox{-minor of}\ V\right\}$.
  If such number does not exist, we set $s = 0$.
  In the case of $s > 0$,
  we have 
  \begin{gather*}
    \gcd \left\{(3 + s) \times (3 + s) \mbox{-minors of}\ V_{m}\right\} = \left\lvert (4+2m) \right\rvert \gcd \left\{(s \times s) \mbox{-minors of}\ V\right\} \ \mbox{and}  \\ 
    \gcd \left\{(3 + s) \times (3 + s) \mbox{-minors of}\ V_{m}\right\} = \left\lvert (4+2n) \right\rvert \gcd \left\{(s \times s) \mbox{-minors of}\ V\right\}.
  \end{gather*}
  If $s=0$,
  \begin{gather*}
    \gcd \left\{3 \times 3 \mbox{-minors of}\ V_{m}\right\} = \left\lvert (4+2m) \right\rvert \mbox{and} \\
    \gcd \left\{3 \times 3 \mbox{-minors of}\ V_{m}\right\} = \left\lvert (4+2n) \right\rvert.
  \end{gather*}
  Since $\left\lvert 4 + 2m\right\rvert \neq \left\lvert 4 + 2n\right\rvert$,
  it follows that $F_{m} \not\cong F_{n}$.
  Therefore,
  the family $\left\{(F_{n}, {F_{n}}')\right\}_{n \in \mathbb{Z}}$ contains infinitely many mutually distinct pairs.  

  (3) We prove that the family satisfies the claims (i)--(iii).

  (i) The handlebody-knots $N(F_{n})$ and $N(F_{n}')$ are equivalent to $H(D_{n})$ and $H(D_{n}')$, respectively.
  The diagrams $D_{n}$ and $D_{n}'$ are related by a finite sequence of Reidemeister moves for handlebody-knots (Fig.~\ref{Fig:R-moves}).
  By Theorem~\ref{Thm:R-moves_for_Hdlebody-knots},
  it follows that $N(F_{n})$ and $N(F_{n}')$ are equivalent as handlebody-knots.

  (ii) Let $V$ be a Seifert matrix of the Seifert surface $F(D)$.
  For each $k \in \mathbb{Z}$,
  since the spatial surfaces $F_{k}$ and $F_{k}'$ have the same Seifert matrix $$V_{k} =
    \begin{pmatrix}
    0 & 0 & 0 & 0 & \mathbf{0} \\
    0 & 0 & 1 & 0 & \mathbf{0} \\
    -1 & 0 & 4 & 0 & \mathbf{0} \\
    -1 & 0 & 0 & 4 + 2k & \mathbf{0} \\
    \mathbf{0} & \mathbf{0} & \mathbf{0} & \mathbf{0} & V
    \end{pmatrix},$$
  Seifert matrices of them are unimodularly congruent.

  (iii) Take the dihedral quandle $R_{3}$ and the cyclic rack $C_{2}$ and 
  consider the product $R = R_{3} \times C_{2}$ of racks $R_{3}$ and $C_{2}$.
  Then,
  we have \begin{equation*}
    \min \left\{k \in \mathbb{Z}_{>0} \mid \mbox{for any}\ (x, a), (y, b) \in R, S_{(y, b)}^{k}(x, a) = (x, a)\right\} = 2
  \end{equation*}
  and
  $X = R \times \mathbb{Z}_{2}$ is an MGR defined as in Example~\ref{Ex:MGRs}.

  Now, we show 
  $\left\lvert \mathrm{Col}_{X}(D_{n})\right\rvert \neq \left\lvert \mathrm{Col}_{X}(D_{n}')\right\rvert$.
  Take a Y-orientation of $D_{n}$ as shown in Fig.~\ref{Fig:R_coloring}, and 
  any coloring of the Y-oriented diagram $D_{n}$ is given by Fig.~\ref{Fig:R_coloring} with 
  the following relations obtained from coloring conditions at the vertices $v_{1}$ and $v_{2}$.    
  \begin{figure}[h]
    \centering
    \includegraphics{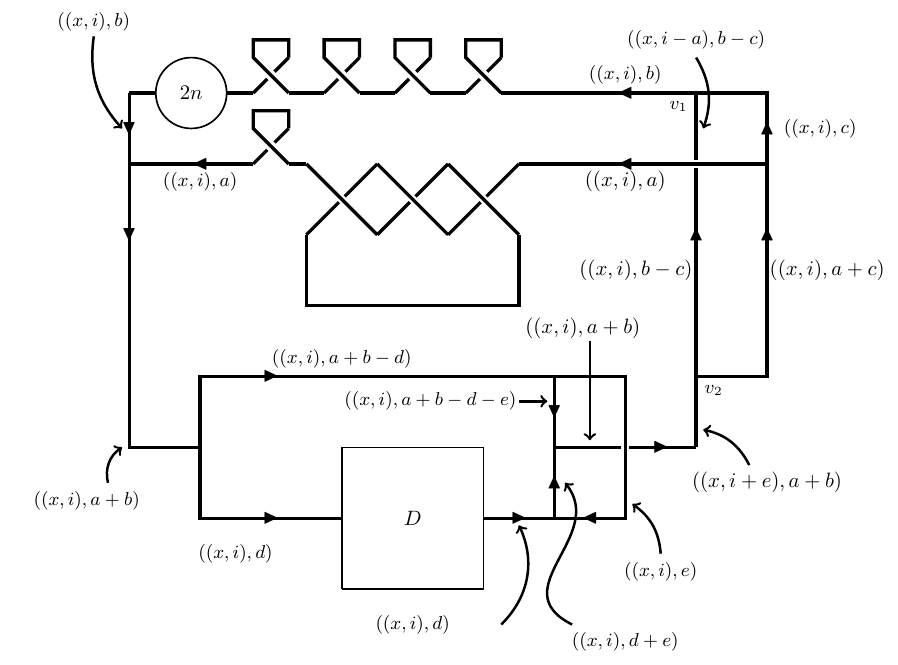}
    \caption{The Y-oriented diagram $D_{n}$ with an $X$-coloring}{\label{Fig:R_coloring}}
  \end{figure}
  \begin{equation*}
    \begin{cases}
      i-a = i & (\mbox{mod}\ 2)\\
      i + e = i & (\mbox{mod}\ 2),
    \end{cases}
  \end{equation*}
  i.e.,
  \begin{equation*}
    \begin{cases}
      a = 0 & (\mbox{mod}\ 2)\\
      e = 0 & (\mbox{mod}\ 2).
    \end{cases}
  \end{equation*}
  For any $((x, i), d) \in X$, we denote the number of $X$-colorings of the sub-diagram with an $X$-coloring illustrated in Fig.~\ref{Fig:Col_D} by $\displaystyle\mathrm{col}_{X}(D; ((x,i),d))$.
  \begin{figure}[h]
    \centering
    \includegraphics{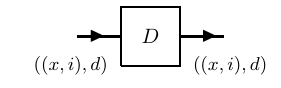}
    \caption{The sub-diagram of $D_{n}$ or $D_{n}'$ with an $X$-coloring}{\label{Fig:Col_D}}
  \end{figure}
  
  Then, we have
  \begin{align*}
    \left\lvert \mathrm{Col}_{X}(D_{n})\right\rvert & = \sum_{((x,i), (b, c, d)) \in R \times \mathbb{Z}_{2}^{3}} \mathrm{col}_{X}(D; ((x, i), d))\\
    & = 4 \left( \bigsqcup_{((x,i), d) \in X} \mathrm{col}_{X}(D; ((x, i), d))\right). 
  \end{align*}
  We remark that $\left\lvert \mathrm{Col}_{X}(D_{n})\right\rvert > 0$ because
  $\mathrm{col}_{X}(D; ((x, i), 0)) > 0$ for any $(x, i) \in R$. 
  Take a Y-orientation of $D_{n}'$ as shown in Fig.~\ref{Fig:col_D2},
  and any coloring of the Y-oriented diagram $D_{n}'$ is given by Fig.~\ref{Fig:col_D2}
  with the following relations obtained from coloring conditions at the vertices $v_{1}$ and $v_{2}$.
  \begin{figure}[h]
    \centering
    \includegraphics{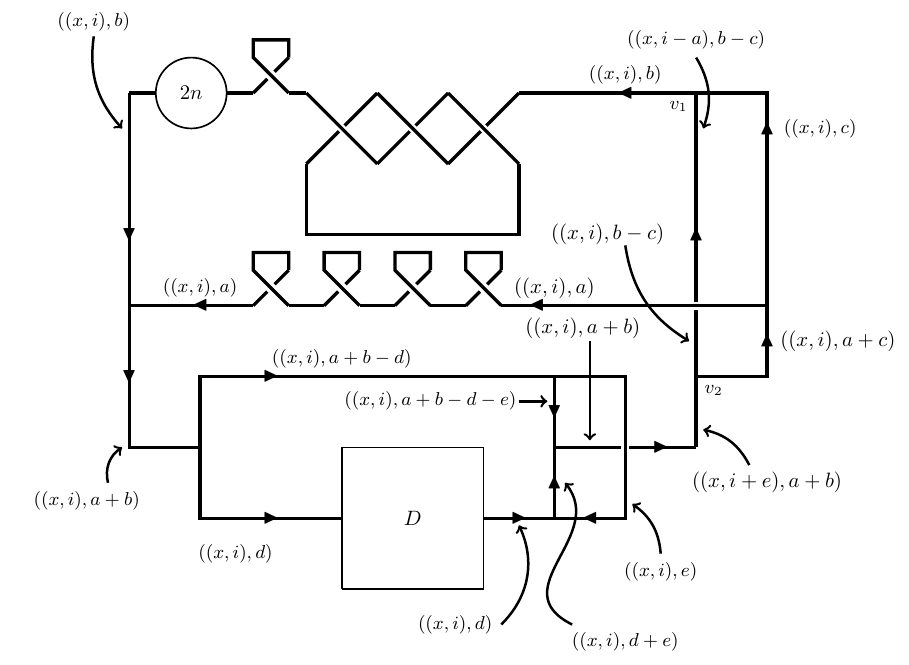}
    \caption{The Y-oriented diagram $D_{n}'$ with an $X$-coloring}{\label{Fig:col_D2}}
  \end{figure}
    \begin{equation*}
    \begin{cases}
      i-a = i & (\mbox{mod}\ 2)\\
      i + e = i & (\mbox{mod}\ 2),
    \end{cases}
  \end{equation*}
  i.e.,
  \begin{equation*}
    \begin{cases}
      a = 0 & (\mbox{mod}\ 2)\\
      e = 0 & (\mbox{mod}\ 2).
    \end{cases}
  \end{equation*}
  When $b = 1$,
  noting that the trefoil knot is $3$-colorable,
  the number of colorings is 
  $\left\lvert \mathrm{Col}_{X}(D_{n}')\right\rvert = 2 \left\lvert \mathrm{Col}_{X}(D_{n})\right\rvert$.
  Therefore, we conclude that $F_{n} \not\cong F_{n}'$.
\end{proof}

\section*{Acknowledgement}

The author would like to thank Keisuke Himeno, Ryoya Kai, Yuko Ozawa, and Yuta Taniguchi for their constant encouragements.
This work was supported by JST SPRING, Grant Number JPMJSP2138.


\vspace{-3mm}

\address{(K. Arai) Department of Mathematics, Graduate School of Science, The University of Osaka, 1-1, Machikaneyama, Toyonaka, Osaka, 560-0043, Japan}

\email{u068111h@ecs.osaka-u.ac.jp}


\end{document}